# Inductive construction of 2-connected graphs for calculating the virial coefficients

E Androulaki[1], S Lambropoulou[2], I G Economou[3] and J H Przytycki[4]

1. National Center for Scientific Research "Demokritos", Institute of Physical Chemistry, GR-15310 Aghia Paraskevi, Greece

2. National Technical University of Athens, Department of Mathematics, Zografou Campus, GR-157 80 Athens, Greece,

3. National Center for Scientific Research "Demokritos", Institute of Physical Chemistry, GR-15310 Aghia Paraskevi, Greece and The Petroleum Institute, Department of Chemical Engineering, PO Box 2533, Abu Dhabi, United Arab Emirates

4. The George Washington University, Department of Mathematics, Washington, DC 20052, USA,

Corresponding author: S Lambropoulou
e-mail: sofia@math.ntua.gr

**Abstract.** In this paper we give a method for constructing systematically all simple 2-connected graphs with $n$ vertices from the set of simple 2-connected graphs with $n$-1 vertices, by means of two operations: subdivision of an edge and addition of a vertex. The motivation of our study comes from the theory of non-ideal gases and, more specifically, from the virial equation of state. It is a known result of Statistical Mechanics that the coefficients in the virial equation of state are sums over labelled 2-connected graphs. These graphs correspond to clusters of particles. Thus, theoretically, the virial coefficients of any order can be calculated by means of 2-connected graphs used in the virial coefficient of the previous order. Our main result gives a method for constructing inductively all simple 2-connected graphs, by induction on the number of vertices. Moreover, the two operations we are using maintain the correspondence between graphs and clusters of particles.

# 1. Introduction

The calculation of thermodynamic properties of non-ideal gases and liquids is usually based on the *equation of state*, which is typically a function of the pressure $p$, the temperature $T$ and the density $\rho$ (or volume $V$) of the system. The general form of an equation of state is $f(p,T,V) = 0$. The most widely used form of an equation of state, used for the calculation of the thermodynamic properties of non-ideal fluids in low or medium pressures, is the *virial equation of state*, which has the general form:

$$\frac{p}{kT} = \rho + B_2(T)\rho^2 + B_3(T)\rho^3 + ...$$

where $B_2(T)$, $B_3(T)$, etc are the 2nd, 3rd, etc *virial coefficients*. Virial coefficients are functions of the temperature for a pure component (and of the composition in the case of mixtures) and can be calculated, in principle, using methods of Statistical Mechanics [McQuarrie D A 1976]. As it turns out, the virial coefficient $B_n$ is an integral of a weighted sum of all 2-connected graphs of $n$ vertices, where particles correspond to vertices and interactions between particles correspond to edges.

A *k-connected graph* is a graph with the property that for any pair of vertices there are at least $k$ disjoint paths connecting them. In order to find all $k$-connected graphs of $n$ vertices one has to extract them from the set of all connected graphs of $n$ vertices. As the number of connected and $k$-connected graphs increases exponentially with $n$, it can be very tedious to isolate the $k$-connected graphs. In the literature there are several results by which we can obtain all $k$-connected graphs by simpler $k$-connected graphs. All these results are basically resting on the operations *addition of an edge* and *amalgamation of vertices*.

$k$-connected graphs are used in various areas apart from Graph Theory and the calculation of virial coefficients in Statistical Mechanics. For example in networks, where a graph is seen as a representation of a network [Bondy J S, Murty U S R 1976] or in Chemistry, where rigid symmetries of 3-connected graphs and, in general, of $k$-connected graphs play an important role [Flapan E 2000].

Motivated by the use of graphs in the calculation of the virial coefficients, we prove in this paper the following result:

**Theorem 1.** *In a simple 2-connected graph of n vertices there is at least one vertex that can be subtracted or removed so that the result will be a simple 2-connected graph with n-1 vertices. Equivalently, every simple 2-connected graph with n vertices arises from a simple 2-connected graph with n-1 vertices by adding a vertex or by subdividing an edge.*

Although the above seems natural we were not able to find anything similar in the literature. Firstly, our two operations, *addition of a vertex* and *subdivision of an edge*, are only related to vertices. Secondly, our method of building 2-connected graphs is inductive with *induction on the number of vertices* and *with induction step equal to* 1. On the contrary, all inductive results in the literature that characterize $k$-connectedness (in general) are basically resting on the operation *addition of an edge* and *the induction step is not necessarily equal to* 1 (see Section 4). Theorem 1 guarantees that, using the two operations addition of a vertex and subdivision of an edge, we obtain the full list of simple 2-connected graphs with $n$ vertices from those with $n-1$ vertices. But, we note that in this list there will be repetitions, that is, pairs of isomorphic graphs (see Remark 3).

The two operations in Theorem 1 correspond naturally to situations in non-ideal fluids. Indeed, looking at a 2-connected graph as a cluster of particles [McQuarrie D A 1976, Mayer J E & Mayer M G 1940], *adding a vertex* means that a particle comes sufficiently close to a cluster of particles and starts to interact with at least two of them, while *subdividing an edge* means that a particle comes close to a cluster in such a way that it pushes away two particles of the cluster, so that they stop interacting with each other, and the new particle begins to interact with these two particles (only). Moreover, we consider only simple graphs, as it has no physical sense to say that two particles interact twice, therefore, we cannot allow more that at most one edge per pair of vertices. Finally, the two operations are particularly adapted to computing the virial coefficients sequentially, since induction on the number of vertices allows stopping the calculation of the virial coefficients on the desired order, the number of vertices, which depends on the properties of the given system.

Our results should be of interest to graph theorists, physical chemists, chemical engineers and electrical engineers. For example, they can be applied to the theory of networks in which, for safety reasons, at least two disjoint paths are required between two nodes (see [Bondy J S, Murty U S R 1976], §3.3). One related work is [Šuvakov M, Tadić B 2006] where a new algorithm is given for growing a planar cellular graph by proper attachment of polygons of length chosen from a given distribution. The resulting category of planar graphs are 2-connected (if "open") or 3-connected (if "closed"), resembling soap froths or nano-particles self-assembled through nonlinear dynamic processes. In the language of our operations, one could start with a triangle and do subdivision of edges, then addition of a vertex of degree 2, followed by possible subdivision of the new edges, and so on. A more applied example is the formation of clusters of particles, met in many physical, chemical or biological situations. In particular, formation of colloids with valency or patchy colloids, which can be used as building blocks of specifically designed self-assembled structures and new colloidal molecules, an extremely new and fast-growing topic. See [Zaccarelli E 2007].

The paper is organized as follows. In Section 2 we explain the theory of virial coefficients and how these are related to 2-connected graphs. In Section 3 we recall some definitions from graph theory and in Section 4 we give known theoretical results on *k*-connected graphs. In Section 5 we describe the operations on 2-connected graphs that we use and we prove that applying these operations yields 2-connected graphs. Finally, we prove our main result (Theorem 1). We extracted our proofs from known theorems in graph theory, which, nevertheless, are very complicated, they are much more general, and they require an extensive knowledge of graph theory. For the extent needed in this paper we give simple and straightforward proofs. Moreover, in the Appendix we give an alternative proof for Theorem 1 (Theorem 2) after extending slightly its statement.

**2. Calculation of the virial coefficients in the classical limit**

The starting point for the calculation of virial coefficients is the *grand-canonical partition function* $\Xi$. The function $\Xi$ can be related to the canonical partition function, $Q$, through the equation:

$$\Xi(V,T,\mu) = \sum_{N=0}^{\infty} Q(N,V,T)\lambda^N \qquad (1)$$

where $N$ is the *number of molecules*, $\mu$ is the *chemical potential*, $\Xi(V,T,\mu)$ is the *grand-canonical partition function*, $Q(N,V,T)$ is the *canonical partition function* and $\lambda = \exp(\beta\mu)$ the *activity*. For $N=0$, the system has only one state with $E=0$, so $Q(N=0,V,T) = 1$. This allows us to write eq. 1 in the following form:

$$\Xi(V,T,\mu) = 1 + \sum_{N=1}^{\infty} Q_N(V,T)\lambda^N \qquad (2)$$

where $Q_N(V,T)$ is equivalent to $Q(N,V,T)$. Using this form of $\Xi(V,T,\mu)$ we can find a suitable form of the *characteristic thermodynamic equation* in the grand-canonical ensemble which can be compared with the virial equation of state, in order to get an expression for the virial coefficients, see [McQuarrie D A 1976], p. 224-226:

$$pV = kT \ln \Xi \tag{3}$$

The procedure described above gives the pressure $p$ as a power series of the density $\rho$:

$$\frac{p}{kT} = \rho + B_2(T)\rho^2 + B_3(T)\rho^3 + ... \tag{4}$$

where

$$B_2(T) = -b_2 = -(2!V)^{-1}(Z_2 - Z_1^2) \tag{5a}$$

$$B_3(T) = 4b_2^2 - 2b_3 = -\frac{1}{3V^2}\left[V(Z_3 - 3Z_2Z_1 + 2Z_1^3) - 3(Z_2 - Z_1^2)^2\right] \tag{5b}$$

...

The $B_i(T)$ $i = 2,3,...$ are the virial coefficients and $Z_i$ are configurational integrals discussed below.

Eqs 5 become significantly more complicated for higher order virial coefficient Nevertheless, the first few virial coefficients are sufficient for the calculation of the equation of state for pressures up to a few hundred atmospheres, as shown in table 1.

| $p$ (atm) | $p/\rho kT$ | | | |
|---|---|---|---|---|
| | $1 + B_2\rho$ | $+ B_3\rho^2$ | $+$ | residual |
| 1 | $1 - 0.00064$ | $+ 0.00000$ | $+...$ | $(+0.00000)$ |
| 10 | $1 - 0.00648$ | $+ 0.00020$ | $+...$ | $(-0.00007)$ |
| 100 | $1 - 0.06754$ | $+ 0.02127$ | $+...$ | $(-0.00036)$ |
| 1000 | $1 - 0.38404$ | $+ 0.68788$ | $+...$ | $(+0.37232)$ |

*Source:* Donald A. McQuarrie, *Statistical Mechanics*, New York: Harper & Row

**Table 1.** The contribution of the first terms in the expansion of $p/\rho kT$ for Argon at 25°C.

This work focuses on the virial coefficients in the thermodynamic limit where the volume becomes arbitrary large, so that our system approximates the properties of a macroscopic system. In the analysis below we will refer only to monoatomic fluids. So, in the classical limit, the canonical partition function $Q_N$ is given by the formula:

$$Q_N = \frac{1}{N!}\left(\frac{2\pi mkT}{h^2}\right)^{3N/2} Z_N \tag{6}$$

Here $h$ is the Planck's constant and $Z_N$ is the *configuration integral*:

$$Z_N = \int ... \int e^{-U_N/kT} d\mathbf{r}_1 d\mathbf{r}_2 ... d\mathbf{r}_N \tag{7}$$

where $r_i$ is the position of the $i$-th molecule. We saw above in eqs. 5 and 6 that for the calculation of any $B_M$ we need the partition functions $Q_N$ for $N \leq M$ or, equivalently, the configuration integral $Z_N$ for $N \leq M$. For example, for the calculation of the second and third virial coefficient we need the configuration integrals:

$$Z_1 = \int d\mathbf{r}_1 = V$$
$$Z_2 = \iint e^{-U_2/kT} d\mathbf{r}_1 d\mathbf{r}_2 \tag{8}$$

and

$$Z_3 = \iiint e^{-U_3/kT} d\mathbf{r}_1 d\mathbf{r}_2 d\mathbf{r}_3$$

Moreover, as we notice, for the calculation of the second virial coefficient we need the potential function $U_2$. For monoatomic molecules, it is logical to assume that the potential function $U_2$ depends only on the distance between the two molecules, so $U_2 = u(r_{12})$, where $r_{12} = |\mathbf{r}_2 - \mathbf{r}_1|$. We then have the expression of $B_2(T)$ as a function of $u(r_{12})$ and by replacing the $Z_1$ and $Z_2$ in eq. 4a:

$$B_2(T) = -\frac{1}{2V}(Z_2 - Z_1^2) = -\frac{1}{2V}\iint \left[ e^{-u(r_{12})/kT} - 1 \right] d\mathbf{r}_1 d\mathbf{r}_2 \tag{9}$$

Using the same approach for the calculation of the third virial coefficient, we need the potential function $U_3(\mathbf{r}_1, \mathbf{r}_2, \mathbf{r}_3)$. The most common approximation for the $U_3$ is:

$$U_3(\mathbf{r}_1, \mathbf{r}_2, \mathbf{r}_3) \approx u(r_{12}) + u(r_{13}) + u(r_{23})$$

that is, we take the potential for three molecules to be the sum of the three pairwise potentials. By replacing $Z_1$, $Z_2$ and $Z_3$ in eq. 4b we obtain the following expression for $B_3(T)$ as a function of $u(r_{12})$, $u(r_{13})$ and $u(r_{23})$:

$$B_3(T) = -\frac{1}{3V^2}\left[ V(Z_3 - 3Z_2 Z_1 + 2Z_1^3) - 3(Z_2 - Z_1^2)^2 \right] =$$
$$= -\frac{1}{3V}\iiint \left[ e^{-u(r_{12})/kT} - 1 \right]\left[ e^{-u(r_{13})/kT} - 1 \right]\left[ e^{-u(r_{23})/kT} - 1 \right] d\mathbf{r}_1 d\mathbf{r}_2 d\mathbf{r}_3 \tag{10}$$

The terms in the square bracket in eqs. 9 and 10 appear very often in equations of the theory of non-ideal gases and are known as the *f-Mayer function*. The *f*-Mayer function is defined as:

$$f_{ij} = f(r_{ij}) = e^{-u(r_{ij})/kT} - 1$$

where $r_{ij} = |\mathbf{r}_j - \mathbf{r}_i|$. Obviously $f(r_{ij}) \to 0$ as $r_{ij} \to \infty$, since $u(r_{ij}) \to 0$ as $r_{ij} \to \infty$. Therefore, eqs. 9 and 10 become:

$$B_2(T) = -\frac{1}{2V} \iint f_{12} d\mathrm{r}_1 d\mathrm{r}_2$$

$$B_3(T) = -\frac{1}{3V} \iiint f_{12} f_{13} f_{23} d\mathrm{r}_1 d\mathrm{r}_2 d\mathrm{r}_3$$

If we look at the integral in $B_3(T)$ more closely we see that it includes three molecules and since $f_{ij} \to 0$ as the *i* and *j* molecules draw away, the product $f_{12} f_{13} f_{23}$ disappears unless the three molecules are close together.

One way to represent schematically the integral in the expression of $B_3(T)$ is the following: we draw a numbered circle for each different index that appears in the product and a line for each pair of molecules that are connected through an *f*-function. The result is a *labeled graph*. For example, the integrals in the second and third virial coefficients can be represented schematically as it is indicated in figure 1a. Now, since each integral represents a cluster of molecules, we call this type of diagrams, like those in figures 1a and 1c, *cluster diagrams of molecules*. The general result of the theory of non-ideal fluids reduces the calculation of virial coefficients to finding cluster diagrams. This very important result gives each virial coefficient through the formulas below:

$$B_{j+1} = \frac{-j}{j+1} \beta_j \qquad \text{(a)} \qquad (11)$$

where

$$\beta_j = \frac{1}{j!V} \int \ldots \int S'_{1,2,\ldots,j+1} d\mathrm{r}_1 d\mathrm{r}_2 \ldots d\mathrm{r}_{j+1} \qquad \text{(b)}$$

and where $S'_{1,2,\ldots,j+1}$ is the sum of all products of *f*-functions that connect the 1, 2, …, *j*+1 molecules, in such a way that the removal of any vertex and the incident edges will not give a disconnected graph. Such graphs are called 2-*connected graphs*. All graphs illustrated in figures 1a and 1c are 2-connected graphs, while those in figure 1b do not have this property. For example, we have $S'_{1,2}$ and $S'_{1,2,3}$ as illustrated in figure 1a respectively, while $S'_{1,2,3,4}$ is the sum of all labeled graphs illustrated in figure 1c. Moreover, if we forget the labelings we may write:

$$S'_{1,2} = -$$
$$S'_{1,2,3} = \triangle \qquad (12)$$
$$S'_{1,2,3,4} = 3\square + 6\boxslash + \boxtimes$$

where the coefficients represent the number of different labelings of the graph. A detailed discussion of the above ideas are given, for example, in [McQuarrie D A 1976].

**Remark 1.** The number of different labelings of an unlabeled graph *G* on *n* vertices is equal to *n*! and inverse-proportional to the order of its automorphism group, Aut(*G*), which detects the symmetries of *G*, see, for example, [BehzadM, Chartrand G, Lesniak-Foster L 1979], Chapter 9, [Godsil C, Royle G 2001], Chapter 2 or [Bollobás B 2001], Chapter 9. Finding Aut(*G*) is not always easy and it is a widely researched area of Graph Theory. Then, according to the type of the potential between the particles, one needs to calculate the value of the integral of S over all the possible positions of the particles. But this is just an operational issue and there are routine algorithms that can do this.

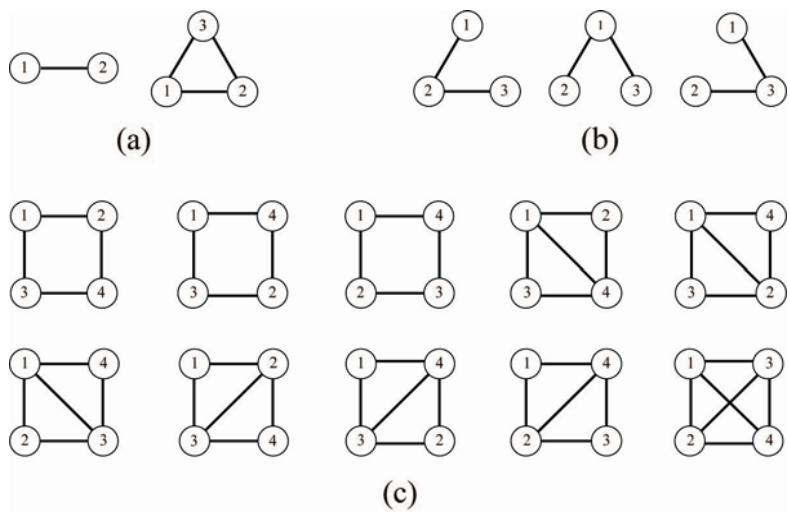

**Figure 1.** Examples of cluster diagrams. (a) The integrals in the second and third virial coefficient. (b) The three topologically equivalent connected diagrams of three molecules. (c) The three different types of 2-connected graphs for four molecules.

To sum up, we have seen that in order to compute virial coefficients we need to know the different types of 2-connected graphs and for each type the number of different labelings. We will explain all these notions in more detail below.

### 3. Definitions from Graph Theory

A *graph G* is defined by the set $V(G)$ of its *vertices*, the set $E(G)$ of its *edges* and a *relation of incidence $f_G$*, which associates each edge with one or two vertices called its *ends*. Equivalently, we can say that it associates each vertex with its adjacent vertices. The *degree* of a vertex is the number of incident edges. A graph $G$ is called *simple* if between two vertices there is at most one edge, and loops are excluded. Throughout this work we will consider only simple graphs. A graph $H$ is said to be a *subgraph* of a graph $G$ if:

$$V(H) \subseteq V(G), \qquad E(H) \subseteq E(G)$$

where $V(G)$ is the *set of the vertices* of a graph $G$ and $E(G)$ is the *set of the edges* of a graph $G$.

In order to find all different 2-connected graphs it is of high importance to find a way to separate graphs. Let $G_1$ and $G_2$ be two graphs with $V(G_1)$ and $V(G_2)$ the corresponding sets of vertices, respectively. The graphs $G_1$ and $G_2$ will be *isomorphic* if there is a bijection $g : V(G_1) \to V(G_2)$ such that two vertices $v_1$, $v_2 \in V(G_1)$ are adjacent in $G_1$ if and only if the vertices $g(v_1)$, $g(v_2) \in V(G_2)$ are adjacent in $G_2$. For example, for three vertices there are two non-isomorphic connected graphs, $\angle$ and $\triangle$, and only one 2-connected graph, namely $\triangle$. Two isomorphic graphs are of the same type.

A *walk* in a graph is a finite non-empty sequence $W = v_0 e_1 v_1 ... e_k v_k$, whose terms are alternately vertices and edges. If the vertices and the edges of a walk $W$ are distinct, $W$ is called a *path*. Two paths of a graph are called *disjoint* if they do not share an inner vertex. A graph is *connected* if for each pair of vertices there is a path joining them. We define an *n-arc* as a graph with *n* edges and (*n*+1) vertices, having the following property: the edges can be enumerated as $A_1, A_2, ..., A_n$, and the vertices $a_0, a_1, a_2, ..., a_n$, in such a way that the ends of $A_j$ are $a_j$ and $a_{j-1}$, for each $j \leq n \in \mathbb{N}$.

The *connectivity* of a graph $G$, $\kappa(G)$, is the smallest number of vertices that must be deleted, together with the incident edges, so that the resulting graph will be disconnected or isomorphic to the isolated vertex. A graph *G* is called *k-connected* if $\kappa(G) \geq k$. Note that *1-connected graph* is the same as *connected graph*. Figure 1a illustrates the only simple 2-connected graphs of two and three vertices. By convention, the graph $K_2$ on two vertices and one edge is 2-connected. Figure 1c illustrates the three different types, up to isomorphism, of simple 2-connected graphs with four vertices, namely, □, ◻ and ⊠. Each of these types has a number of possible labelings, where by a *labelling* is understood a numbering of the vertices. For example, the unlabeled graph □ has three different labellings, as shown in fig. 1c. But how can we distinguish the different labelings of a graph? Let $G_1$ and $G_2$ be two labeled graphs with $f_{G_1}$ and $f_{G_2}$ the corresponding relations of incidence. We say that the graphs $G_1$ and $G_2$ have *different labelings* if they are isomorphic as unlabelled graphs but $f_{G_1} \neq f_{G_2}$ as labeled graphs. In figure 1c we can see the three different types of 2-connected graphs with four vertices.

**4. Results on *k*-connected graphs**

Substantial work has been reported in the literature on *k*-connected graphs. A keystone result in the theory is the following [Menger K 1927]:

**Menger's Theorem.** *A graph G with number of vertices* $n \geq k+1$ *is k-connected if and only if any 2 distinct vertices of G are connected by at least k disjoint paths.*

It follows directly that each vertex of a 2-connected graph is of degree at least 2.
Further, quoting Bollobás from [Bollobás B 1978]: "*One of the most important goals of the theory of k-connected graphs is to compile a list of all k-connected graphs. A natural way of achieving this would be to give some operations producing k-connected graphs from k-connected graphs such that every k-connected graph can be obtained from certain simple k-connected graphs by repeated applications of the operations.*". Indeed, this has been accomplished by Tutte [Tutte W T 1961] for 3-connected graphs and by Slater [Slater P J 1974] for 4-connected graph. In the case of 4-connected graphs, Slater found rather complicated and not too informative operations. On the other hand, Tutte proved the following result: "*A graph is 3-connected if and only if it is a wheel or can be obtained from a wheel by repeated applications of the following two operations: a) The addition of an edge, b) The replacement of a vertex x of degree* $\geq$ *4 by two adjacent vertices x', x'' and joining every neighbor of x to exactly one of x', x'' in such a way that both x' and x'' will have degree* $\geq$ *3. (This operation is called splitting of the vertex x.)*"

Moreover, it has been common sense [Bollobás B 1978, Plummer M D 1968, Dirac G A 1967] that since the addition of an edge does not decrease the connectivity of a graph, in order to describe all *k*-connected graphs, it suffices to describe all *minimally k-connected graphs*, that is, graphs that are *k*-connected but loose this property if we remove any of their edges. The strongest results are by Halin and Mader, who worked on minimally *k*-connected graphs. The first results were obtained by Halin [Halin R 1969] for example: "*Every finite k-connected graph G contains either a vertex of valency k or an edge e such that the graph arising from G by the deletion of e remains k-connected*". Halin also conjectured a number of extensions of those results. Subsequently all conjectures were proved by Mader [Mader W 1971, Mader W 1978]. Halin's result is as follows.

We will focus here on 2-connected graphs. Substantial work has been reported in the literature on 2-connected graphs, see for example [Bondy J S, Murty U S R 1976], [Diestel R 2000], [Harary F 1969] and [Tutte W T 2001], Chapter III. Dirac [Dirac G A 1967] and Plummer [Plummer M D 1968] have proved the following inductive result: *"For each i, $0 \leq i \leq k$ ($k \geq 1$) let $G_i$ be an edge $x_i x_{i+1}'$ or a minimally 2-connected graph containing compatible vertices $x_i$ and $x_{i+1}'$. Let G be the graph obtained from $\bigcup_0^k G_i$ by identifying $x_i$ with $x_{i+1}'$ for $0 \leq i \leq k$ and joining $x_0$ to $x_{k+1}'$. Then G is minimally 2-connected. Conversely, every k-connected graph can be obtained in the way described above."* As we see, they obtain *all minimally 2-connectetd graphs from minimally 2-connected graphs*. Also, the operation used is *amalgamation* of minimally 2-connected graphs and edges (via vertices). Another inductive result about construction of 2-connected graphs is by Diestel [Diestel R 2000]: *"A graph is 2-connected if and only if it can be constructed from a cycle by successively adding arcs, which have their endpoints in the already constructed graph"*, for a proof see, for example, [Diestel R 2000], p. 44. Diestel uses the operation of *adding arcs*, that is, connected graphs with two endpoints and all the other vertices of degree 2.

Finally, we have the following results, which form the basis of our work.

**Tutte's Theorem.** *Let G be a simple 2-connected graph having at least two edges. Then G can be represented as a union of a simple 2-connected subgraph H and an arc L that avoids H but has both its ends in H.* For a proof of Tutte's Theorem see, for example, [Tutte W T 2001], p.57.

**Corollary 1.** *Let G be a simple 2-connected graph. Then either there is a vertex in G of degree 2 or there is an edge e such that $G - e$ is 2-connected.*

**5. Inductive derivation of 2-connected graphs based on vertices**

The main problem in finding all 2-connected graphs with *n* vertices is to separate them from the connected graphs with *n* vertices. It is obvious that the set of the 2-connected graphs with *n* vertices is a genuine subset of all connected graphs with *n* vertices. In table 2 we can see the number of connected graphs with *n* vertices $C(n)$ and the corresponding number of 2-connected graphs with *n* vertices $S(n)$, for $n \leq 19$. For a full list of 2-connected graphs with $n \leq 7$, see, for example, the Appendix of [Hoover W G, De Rocco A G 1961].

| $C(n)$ | $n$ | $S(n)$ |
|---|---|---|
| 1 | 1 | 0 |
| 1 | 2 | 1 |
| 2 | 3 | 1 |
| 6 | 4 | 3 |
| 21 | 5 | 10 |
| 112 | 6 | 56 |
| 853 | 7 | 468 |
| 11117 | 8 | 7123 |
| 261080 | 9 | 194066 |
| 11716571 | 10 | 9743542 |
| 1006700565 | 11 | 900969091 |
| 164059830476 | 12 | 153620333545 |
| 50335907869219 | 13 | 48432939150704 |
| 29003487462848061 | 14 | 28361824488394169 |
| 31397381142761241960 | 15 | 30995890806033380784 |
| 63969560113225176176277 | 16 | 63501635429109597504951 |
| 245871831682084026519528568 | 17 | 244852079292073376010411280 |
| 1787331725248899088890200576580 | 18 | 1783160594069429925952824734641 |
| 24636021429399867655322650759681644 | 19 | 24603887051350945867492816663958981 |

<u>*Source:*</u> *Wolfram MathWorld* (http://mathworld.wolfram.com/), *The On-Line Encyclopedia of Integer Sequences* (http://www.research.att.com/~njas/sequences/)

**Table 2.** The number $C(n)$ of connected graphs and the number $S(n)$ of 2-connected graphs for vertices $n \leq 19$.

We can easily see that these numbers grow exponentially with *n*. Therefore, it becomes extremely difficult to extract the 2-connected graphs from the set of connected graphs. For that, it would be very nice if we could build inductively all 2-connected graphs with *n* vertices using the 2-connected graphs with *n*-1 vertices. With this aim we allow on graphs the following operations.

**Definition 1.** *Addition of a vertex* in a graph *G* is the operation that adds a new vertex *v* of degree at least 2 in the set *V*(*G*) and its incident edges in the set *E*(*G*). We denote this operation as "$G+(v)$". *Subtraction of a vertex* from a graph *G* is the inverse operation, where we remove a vertex *v* with all its incident edges. We denote this operation as "$G-(v)$".

**Definition 2.** *Addition of an edge* in a graph *G* is the operation where we connect two vertices of *G* by a new edge *e*. So, we add *e* in *E*(*G*). We denote this operation as "$G+e$". The inverse operation is called *removal of an edge*. We denote this operation "$G-e$".

**Definition 3.** *Subdivision of an edge* of a graph *G* is the operation where in an edge *e* with ends the vertices $v_1$ and $v_2$ we add a new vertex *v*. So, we have addition of the vertex *v* in the set *V*(*G*) and replacement of the edge *e* from the edges ($v_1$ *v*) and (*v* $v_2$) in the set *E*(*G*). The inverse operation is called *removal of a vertex*. Here, the two edges incident to a vertex *v* of degree 2 are replaced by one edge and the vertex *v* is removed. We denote this operation "$G-v$".

**Remark 2.** It is worth noticing that any two of the above operations can produce the third one. Indeed, addition of a vertex can be seen as the combination of subdividing an edge and adding edges. Whilst, subdivision of an edge can be seen as addition of a vertex of degree 2 combined with the removal of the old edge. Finally, addition of an edge can be seen as addition of vertex of degree 2 and removal of this new vertex.

In this work we focus on the two operations that add a new vertex in an existing graph, that is, "addition of a vertex" and "subdivision with an edge", because our main result (Theorem 1) is stated by means of induction on the number of vertices. Yet, the operation "addition of an edge" helps in the proof of our statements.

We shall now show that the two operations that add a new vertex, applied on a 2-connected graph yield a 2-connected graph with one more vertex.

**Proposition 1.** *In a simple 2-connected graph G with n-1 vertices the subdivision of an edge by a vertex v yields a simple 2-connected graph $G'$ with n vertices.*

*Proof.* Let $G' = G + (v)$ and (ab) be the subdivided edge. We need to show that $\kappa(G') \geq 2$, that is, the subtraction of any vertex $u \in V(G')$ will give a connected graph. Or, equivalently, that between any two vertices in $G'$ there are at least two disjoint paths connecting them. For a pair of vertices $w, w' \in V(G)$ there are at least two disjoint paths in G connecting them, since $G$ is 2-connected. If one of these paths involves the edge (ab) then in $G'$ it will involve instead $(av) \cup (vb)$. So, in $G'$ there are also at least two disjoint paths between $w$ and $w'$.

Consider now the new vertex *v* paired with a vertex *w*. Then, what we know is that between *w* and say *a* there are at least 2 disjoint paths in $G'$, say $\gamma_1$ and $\gamma_2$. If one of them contains the vertex *v*, say $\gamma_1$, then we can consider the paths $\gamma_2 \cup (av)$ and $\gamma_1 \setminus (av)$ which are disjoint and connect the vertices *v* and *w*. If now neither $\gamma_1$ nor $\gamma_2$ contains the vertex *v*, then, obviously, neither contains the vertex *b* as well. We can also find two disjoint paths connecting *v* and *w*. This is true because if we look at the vertices *b* and *w* in *G*, we see that they connect through the path $\gamma_1 \cup (ab)$ and as *G* is 2-connected, there is at least one more disjoint path, say $\gamma_3$, that connects them as well. Obviously, $\gamma_1$ and $\gamma_3$ are disjoint. So, we can find in $G'$ at two disjoint paths connecting *v* and *w*, namely $(va) \cup \gamma_1$ and $(vb) \cup \gamma_3$. Thus, $G'$ is 2-connected. □

**Proposition 2.** *In a simple 2-connected graph G with n-1 vertices the addition of a vertex yields a simple 2-connected graph with n vertices.*

*Proof.* Let *v* be a new vertex and let $G' = G + (v)$. We need to show that $\kappa(G') \geq 2$, that is, the subtraction of any vertex will give a connected graph. We have the following two cases. First, suppose that from the graph $G'$ is subtracted the vertex *v*. The remaining graph is *G*, which is 2-connected, therefore connected.

Suppose now that from the graph $G'$ is subtracted one of the rest *n*-1 vertices, say *u*. Then the graph that arises is the same as the one that would arise from *G* if we subtracted the vertex *u*, but with the extra vertex *v* of degree at least 2. So, every pair $w, w' \in V(G)$ of old vertices is path-connected. Consider now the vertex *v* paired with a vertex *w*. Then *v* is connected with at least one vertex $w' \in V(G)$ (in fact, *v* connects with at least two, if the vertex *u* is not adjacent to *v*). However, $w'$ is path-connected to *w*, hence *v* is connected to *w*. Therefore, for each pair of vertices there is at least one path that connects them. Thus, by Menger's Theorem $G'$ connected.

Suppose now that from the graph $G'$ is subtracted one of the rest *n*-1 vertices witch is adjacent to the vertex that we have added. Then the resulting graph is the same as the one that would arose from *G* if we would have subtracted the same vertex and before the addition of the new vertex, but with one extra vertex of degree at least 1. So, the new vertex and any other vertex of the new graph, even if it were adjacent to the vertex that we have subtracted, is connected with at least one more vertex. Therefore, for each pair of vertices there is at least one path that connects them. Thus, by Menger's Theorem, $G'$ is connected. □

Note that, in the case where *v* is of degree exactly 2 Proposition 2 follows from Tutte's Theorem.

The question now is whether, using the above results, all 2-connected graphs with *n* vertices arise from the 2-connected graphs with *n*-1 vertices. Indeed we claim the following, which is our main result.

**Theorem 1.** *In a simple 2-connected graph of n vertices there is at least one vertex that can be subtracted or removed so that the result will be a simple 2-connected graph with n-1 vertices. Equivalently, every simple 2-connected graph with n vertices arises from a simple 2-connected graph with n-1 vertices by adding a vertex or by subdividing an edge.*

Using Theorem 1 we have constructed all 2-connected graphs exhaustively up to $n = 7$. For a list of these graphs see Appendix 1 of Chapter 3 of [Androulaki E].

Concerning the operation "removal of an edge", we give the following result which is known (compare with [Tutte W T 2001], p.67, Theorem III.30), but for which we give a more direct proof.

**Proposition 3**. *If from a simple 2-connected simple graph G we remove an edge e the remaining graph $G'$ is either 2-connected or connected. Moreover, if $G'$ is 2-connected, then G must be of the form*:

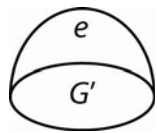

*Iif $G'$ is connected then G must be of the form*:

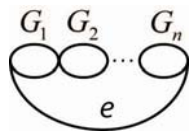

where $G_i$, $i \leq n$ are 2-connected subgraphs with $n \geq 2$.

*Proof.* Let $e \in E(G)$ and $G' = G - e$. We will show that $G'$ is connected. Let $w, w' \in V(G') = V(G)$. By the fact that $G$ is 2-connected there are at least two disjoint paths between $w$ and $w'$, and at most one of them will involve the edge *e*. So, for every pair of vertices in $G'$ there is at least one path connecting them $\underset{\text{Menger}}{\Rightarrow}$ $G'$ is (at least) connected. Obviously, if $G' = G - e$ is 2-connected then $G$ must be isomorphic to the graph:

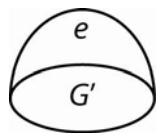

If $G' = G - e$ is not 2-connected then we claim that between the vertices $v_1, v_2$ of the edge $e$ there are no two disjoint paths in $G' = G - e$. Indeed, suppose there are two disjoint paths $\gamma_1, \gamma_2$ connecting $v_1, v_2$. Then we will show that for any pair of vertices $w_1, w_2$ there are two disjoint paths in $G'$, contradiction.

Note first that if between $w_1, w_2$ there are two disjoint paths not involving the edge $e$ there is nothing to show. Suppose now that between $w_1, w_2$ there are exactly two disjoint paths $\delta_1, \delta_2$, one of them, say $\delta_1$, containing the edge $e$. Then $\delta_1$ also contains the vertices $v_1, v_2$. If the path $\delta_2$ involves parts only of, say, $\gamma_2$, then in $\delta_1$ we substitute $e$ by $\gamma_1$ and we are done. If now $\delta_2$ involves parts of both $\gamma_1$ and $\gamma_2$, we do the following:

Let $a$ and $b$ be the first, respectively last intersection of $\delta_2$ with $\gamma_1 \cup \gamma_2$. We distinguish two cases:

If $a, b$ both belong to, say, $\gamma_2$ we define new paths $\zeta_1, \zeta_2$ between $w_1, w_2$ as follows (view figure 2a):

$$\zeta_1 \equiv \left(part\ (w_1, v_1) \in \delta_1\right) \cup \gamma_1 \cup \left(part\ (v_2, w_2) \in \delta_1\right)$$
$$\zeta_2 \equiv \left(part\ (w_1, a) \in \delta_2\right) \cup \left((a, b) \in \gamma_2\right) \cup \left(part\ (b, w_2) \in \delta_2\right)$$

The paths $\zeta_1, \zeta_2$ are obviously disjoint.

Finally, if $a$ belongs to, say, $\gamma_2$ and $b$ belongs to $\gamma_1$ then we define new paths $\zeta_1, \zeta_2$ between $w_1, w_2$ as follows (view figure 2b):

$$\zeta_1 \equiv \left(part\ (w_1, v_1) \in \delta_1\right) \cup \left(part\ (v_1, b) \in \gamma_1\right) \cup \left(part\ (b, w_2) \in \delta_2\right)$$
$$\zeta_2 \equiv \left(part\ (w_1, a) \in \delta_2\right) \cup \left(part\ (a, v_2) \in \gamma_2\right) \cup \left(part\ (v_2, w_2) \in \delta_1\right)$$

The paths $\zeta_1, \zeta_2$ are obviously disjoint.

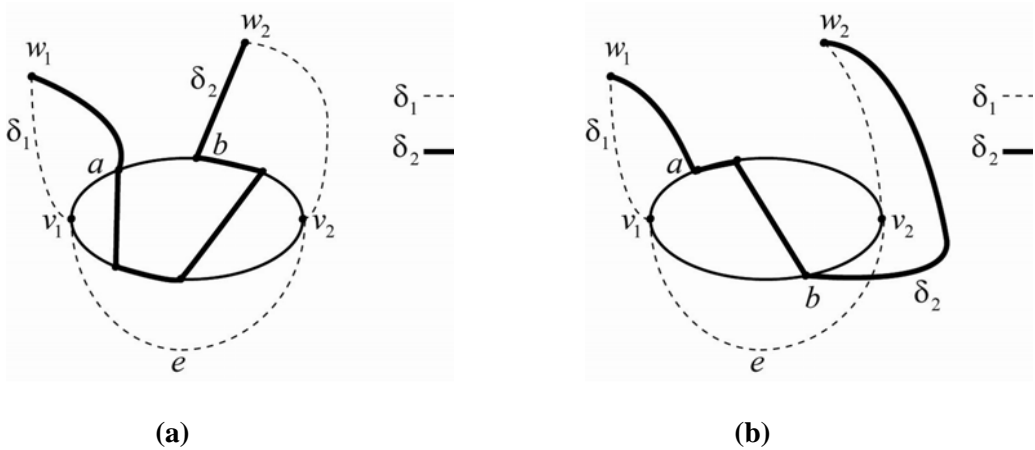

**Figure 2**: Illustration of the paths used in the proof of Proposition 2.

We showed that in any case $G'$ is 2-connected, contradiction. Therefore, $G$ must be of the form:

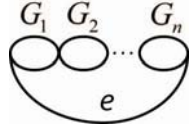

where $G_i$, $i \leq n$ are 2-connected subgraphs with $n \geq 2$. □

Note that, the first case of Proposition 3 is a special case of Tutte's Theorem (where the arc $L$ is only one edge).

We shall finally give two results, which will rest the proof of Theorem 1.

**Proposition 4.** *Let G be a simple 2-connected graph. If there is in G a vertex v of degree 2, then either v can be removed and $G - v$ is 2-connected or v can be subtracted and $G - (v)$ is 2-connected.*

*Proof.* There are two cases for the position of the vertex *v*. Indeed, suppose that only *a* is of degree 2. Then the subtraction of the vertex *b* would give a non-connected graph, which is a contradiction, since *G* is 2-connected. (see figure 3b).

<u>Case 1:</u> The vertices *a* and *b* that are adjacent to the vertex *v*, are connected with an edge (see figure 3a). If the vertices *a* and *b* are of degree 2 then *G* is the triangle *vab* and the statement holds. If the degrees of the vertices *a* and *b* are greater than 2, then there are vertices *c* and *d*, adjacent to the vertices *a* and *b*, where generally $c \neq d$, (see figure 3c). We will show that, in this case, the graph $G - (v)$ is 2-connected. By Menger's theorem we can equivalently show that between any two vertices of the graph there are at least two disjoint paths.

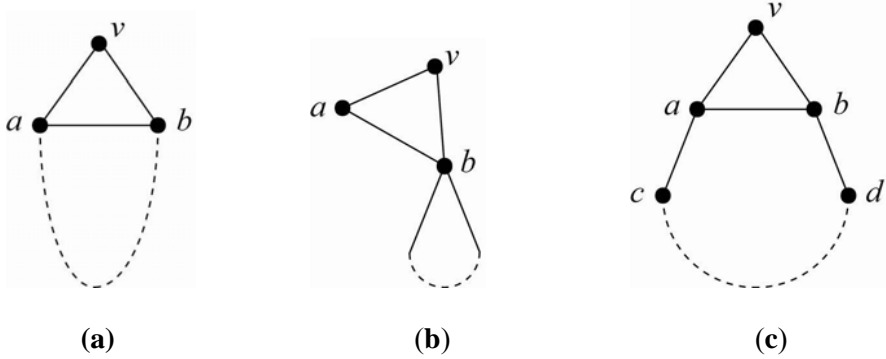

(a)　　　　　　　　(b)　　　　　　　　(c)
**Figure 3**

Let's consider what happens in *G* after the subtraction of the vertex *v*. Take first the vertices *a* and *b*. Between the vertices *a* and *b* there is one path, the edge (*ab*). We must reassure that there is one more path that connects them and is disjoint from the edge (*ab*). To find this path we look at the vertices *c* and *d*. Between them there is the path (*cabd*) in *G* and, according to Menger's Theorem, there is one more path *γ* in *G* disjoint from (*cabd*) that connects the vertices *c* and *d*. Clearly *γ* is in $G - (v)$. Hence, we have the two disjoint paths (*ab*) and $((ac) \cup \gamma \cup (db))$ in $G - (v)$ between *a* and *b*. Take now two vertices $v_1, v_2$ in $G - (v)$, such that it is not the case $v_1 = a, v_2 = b$. Since *G* is 2-connected, there are two disjoint paths between $v_1, v_2$ in *G*. If the one of the two paths contains the vertex *v*, then, it must contain the edges (*av*) and (*vb*). So, the other part cannot contain the edge (*ab*). Therefore, we can replace in the first path the edges (*av*) and (*vb*) by (*ab*), obtaining two disjoint paths in $G - (v)$. Hence, $G - (v)$ is 2-connected.

<u>Case 2:</u> The vertices *a* and *b*, that are adjacent to the vertex *v*, are not connected with an edge. Obviously, the removal of the vertex *v* does not affect the validity of Menger's Theorem in the graph $G-v$ since in each path that contains the edges $(av)$ and $(vb)$, they are just replaced by the new edge $(ab)$. Therefore, in this case also $G-v$ is 2-connected. □

**Proposition 5.** *Let G be a simple 2-connected graph with no vertex of degree 2. Then there is a vertex v in G such that $G-(v)$ is 2-connected.*

*Proof.* We will show this by induction on the number of edges. Since we are dealing with simple graphs, the only graph with no vertex of degree 2 that has the least number of edges is $K_4$, which has six edges:

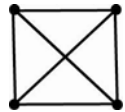

As we can see, in $K_4$ the subtraction of any vertex gives $K_3$, which is a 2-connected graph.

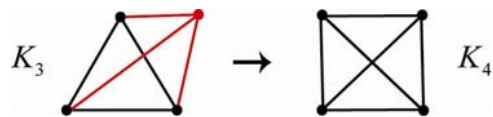

Now, let $n \in \mathbb{N}$, $n > 6$ and suppose that the statement holds for 2-connected graphs with number of edges $E < n$. We will show that the statement holds for 2-connected graphs with $E = n$. Let $G$ be a 2-connected graph with $E(G) = n$ such that every vertex of $G$ is of degree at least 3. Then, from Corollary 1, there is an edge $e = (ab)$ in $G$ such that $G-e$ is 2-connected. We have $E(G-e) = n-1$, and there are two cases for the graph $G-e$.

<u>Case 1:</u> The graph $G-e$ has no vertex of degree 2. Then, by the induction hypothesis, there is a vertex $v$ in $G-e$, such that $(G-e)-(v)$ is 2-connected. If $v \neq a$ or $b$ then $G-(v)$ is obviously also 2-connected. If $v = a$ or $b$ then $(G-e)-(v) = G-(v)$, so $G-(v)$ 2-connected.

<u>Case 2:</u> The graph $G-e$ has at least one vertex $v$ of degree 2 ($v$ must be $a$ or $b$). Then Proposition 4 reassures that $v$ can either be removed and $G-e-v$ is 2-connected or $v$ can be subtracted and $G-e-(v)$ is 2-connected. Therefore, $G-(v)$ is 2-connected. □

We proceed now with the proof of Theorem 1.

***Proof of Theorem* 1.** Let $G$ be a simple 2-connected graph with $n$ vertices. Suppose first that there is a vertex $v$ in $G$ of degree exactly 2. Then, by Proposition 4, $v$ can either be removed and $G-v$ is 2-connected, or $v$ can be subtracted and $G-(v)$ is 2-connected. Suppose now that all vertices in $G$ are of degree greater than 2. Then, by Proposition 5, there is a vertex $v$ that can be subtracted, such that $G-(v)$ is 2-connected. Therefore, every simple 2-connected graph with $n$ vertices arises from a simple 2-connected graph with $n-1$ vertices, by means of either the addition of a new vertex or the subdivision of an edge. □

In [Androulaki E] we demonstrate the application of Theorem 1 for obtaining the full list of 2-connected graphs, with no repetitions, for $n \leq 7$.

**Remark 3.** Theorem 1 guarantees that, using the two operations *addition of a vertex* and *subdivision of an edge*, we obtain the full list of 2-connected graphs with $n$ vertices from those with $n-1$ vertices. But in this list there may be repetitions, that is, pairs of isomorphic graphs. An example is illustrated in figure 4.

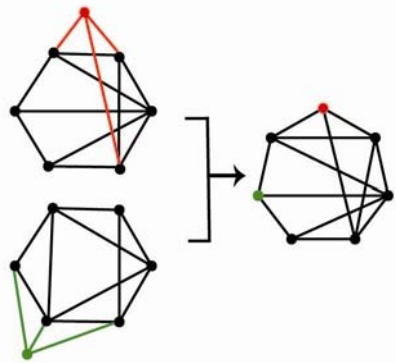

**Figure 4:** Two non-isomorphic 2-connected graphs on six vertices give rise to the same 2-connected graph on seven vertices.

In order to give the full list of 2-connected graphs with no repetitions one has to detect the graphs that are isomorphic and choose a representative for each isomorphism class. This is easy to do by hand for *n* up to 6, see table 3. We have also done it by hand for $n = 7$, see [Androulaki E]. For $n > 7$ one can use invariants of graphs, such as the chromatic polynomial, the Tutte polynomial, etc, see for example [Tutte W T 2001], Chapter IX. Of course, there are algorithms and programs analyzing isomorphisms of graphs. For example, the package NAUTY [McKay B]. None of them, however, has a polynomial complexity. It is an open (and important) problem whether the graph isomorphism problem has a polynomial complexity algorithm, see, for example, [Garey M R, Johnson D S 1979]. This problem is particularly attractive and actively researched by people working on quantum computing: is there a polynomial complexity, quantum algorithm for the graph isomorphism problem?

We believe that our inductive results can be generalized to *k*-connected graphs, with proper adaptation of the operations.

**Conclusions**

- In the present work we have considered only simple graphs. This is because we wanted to preserve the equivalence between graphs and clusters of particles. Indeed, as we have already stated, vertices have a natural representation as particles, whereas edges represent the existence of interaction between two particles. It has no sense to say that two particles interact twice, therefore, we cannot allow more than at most one edge per pair of vertices.

- In order to maintain the naturality of the correspondence between particles and vertices, we allowed operations on graphs that correspond to changes in the state of the clusters. More precisely, the addition of a vertex can be seen as a new particle coming close enough to the cluster in order to interact with at least two existing particles in the cluster. Moreover, subdividing an edge means that a particle comes close enough to two other particles in the cluster and starts interacting with them, pushing them at the same time away, so that they loose their interaction. Finally, addition of an edge means that two particles in the cluster come close enough so that they start to interact. Analogous situations can be considered for the inverse operations.

- According to our results, one can obtain the whole list of simple 2-connected graphs inductively with induction on the number of vertices and with induction step equal to 1, using only the operations in Theorem 1. We extracted our proofs from known theorems in graph theory, which nevertheless are very complicated, they are much more general and they require an extensive knowledge on graphs. For the extent needed in this paper we give much simpler straightforward proofs.

- In the obtained list of simple 2-connected graphs there will be repetitions. To eliminate them, one must detect all isomorphic graphs and then choose one representative for each isomorphism class (see Remark 3). The different 2-connected graphs of up to $n$ vertices are the summands of $S'$ that appear in eq. 11 of the virial coefficients $B_n$. Further, in order to find their coefficients in the combinatorial expression of each $S'$ (recall eqs. 12) one has to find the number of all different labelings for each 2-connected graph (see Remark 1). Then, according to the type of the potential between the particles, one needs to calculate the value of the integral of $S'$ over all the possible positions of the particles. But this is just an operational issue and there are routine algorithms that can do this.

**Acknowledgments.** The authors are pleased to acknowledge a discussion with Béla Bollobás. The last author was partially supported by the George Washington CCAS/UFF summer award, by the GWU REF award and by the NTUA.

**Appendix**

In this Appendix we give an alternative proof of Theorem 1 after extending slightly its statement. We begin with some definitions that will be needed in the sequel.

Let $H$ be a subgraph of a graph $G$. We say that $H$ is an *edge-proper subgraph of $G$* if $E(G) - E(H)$ is non-empty. Similarly, we say that $H$ is a *vertex-proper subgraph of $G$* if $V(G) - V(H)$ is non-empty. From the finiteness of $G$ we have the following:

**Lemma 1:** *Let $G$ be a 2-connected graph.*
   (i)    *If $H$ is a 2-connected vertex-proper subgraph of $G$ then there is a maximal 2-connected vertex-proper subgraph of $G$ (denoted by $H^{v-\max}$) containing $H$.*
   (ii)   *If $H$ is a 2-connected edge-proper subgraph of $G$ then there is a maximal 2-connected edge-proper subgraph of $G$ (denoted by $H^{e-\max}$) containing $H$.*

We will also use the following elementary lemma:

**Lemma 2:** *Let $H$ be a vertex-proper subgraph of a 2-connected graph $G$. Assume also that $H$ has at least two vertices. Then there is a vertex $v$ in $V(G) - V(H)$ and two disjoint paths $\gamma_1$ and $\gamma_2$ connecting $v$ to $H$. "Disjoint" means that $\gamma_1 \cap \gamma_2 = \{v\}$ and "connecting to $H$" means that $\gamma_i \cap H$, $i = 1, 2$ is one vertex.*

Lemma 2 is a special case of Menger's Theorem, but for completeness we give a short proof.

*Proof.* Let $v'$ be a vertex in $V(G) - V(H)$ and connect it to a vertex $v_1$ of $H$ by a path $\gamma_1'$ such that $\gamma_1' \cap H = \{v_1\}$. As $G$ is 2-connected, thus $G - v_1$ is connected. Let $v_2$ be a vertex of $G - (v_1)$ connected to $v'$ in $G - (v_1)$ by a path $\gamma_2'$ such that $\gamma_2' \cap H = \{v_2\}$. The paths $\gamma_1'$ and $\gamma_2'$ are not necessarily disjoint but $\gamma_1' \cap \gamma_2' \subset G - H$. Let $v$ be the last vertex of $\gamma_1'$, when travelling from $v'$, which also belongs to $\gamma_2'$. Obviously $v \in V(G - H)$. Let also $\gamma_i$, $i$=1, 2, be the part of $\gamma_i'$ connecting $v$ to $v_i$. By construction, $\gamma_1 \cap \gamma_2 = \{v\}$ and Lemma 2 is proven. □

We are now ready to give another proof of Theorem 1. In fact we will give a small generalization of it for $n \geq 4$ as follows:

**Theorem 2 (Generalization of Theorem 1).** *Let $G$ be a simple 2-connected graph of at least four vertices. Then either there exists a vertex $v$ such that $G - (v)$ is 2-connected or $G$ has at least four vertices of degree 2.*

*Proof of Theorem 2.* Assume that $G$ contains no vertex such that $G - (v)$ is 2-connected. $G$ is not a tree, so let $C_k$, $k \geq 3$ be the shortest cycle in $G$. If $G = C_k$ then $k \geq 4$ and Theorem is proven. Otherwise, by Lemma 1 there is a maximal 2-connected vertex-proper subgraph $C_k^{v-\max}$ of $G$. By Lemma 2 there is a vertex $v$ in $G - C_k^{v-\max}$ connected to $C_k^{v-\max}$ by two disjoint paths $\gamma_1$ and $\gamma_2$. The graph $C_k^{v-\max} \cup \gamma_1 \cup \gamma_2$ is 2-connected, thus by the maximality of $C_k^{v-\max}$ it contains all vertices of $G$. Furthermore, every vertex of $(\gamma_1 \cup \gamma_2) - C_k^{v-\max}$ is of degree 2. If there were only one such vertex, say $v$, then $G - (v) = C_k^{v-\max}$ is 2-connected, contradiction. Thus $\gamma_1 \cup \gamma_2$ contains at least two vertices of degree 2. Our goal, now, is to show that, if for every vertex *v* the graph $G-(v)$ is not 2-connected, then *G* has at least four vertices of degree 2. To show this we start from the first part of the proof and we consider a path $\gamma_3$ in $C_k^{v-\max}$ that connects the endpoints of $\gamma_1 \cup \gamma_2$. We then have the cycle $C_s = \gamma_1 \cup \gamma_2 \cup \gamma_3$, with $s \geq 4$. As *G* is not a cycle, we can repeat the previous consideration starting from $C_s$ and produce two new vertices of degree 2 outside $C_s$. This completes our proof of Theorem 1. To see that Theorem 2 is a generalization of Theorem 1 consider a simple 2-connected graph with no vertex *v* such that $G-(v)$ is 2-connected. Therefore, as we have shown above, the path $\gamma_1 \cup \gamma_2$ contains at least two vertices of degree 2, so by Proposition 4 one of these vertices can be removed, and we obtain a simple 2-connected graph with one less vertex. □

**Remark 4.** Our generalization of Theorem 1 is the best possible in the sense that there are 2-connected graphs with $G-(v)$ not 2-connected for any vertex *v*, with arbitrarily many vertices and only four vertices of degree 2. See figure 5, where an infinite family of 2-connected graphs is illustrated, with all vertices inessential and four vertices of degree 2. *Inessential vertex* means that its deletion results in a graph which is not 2-connected.

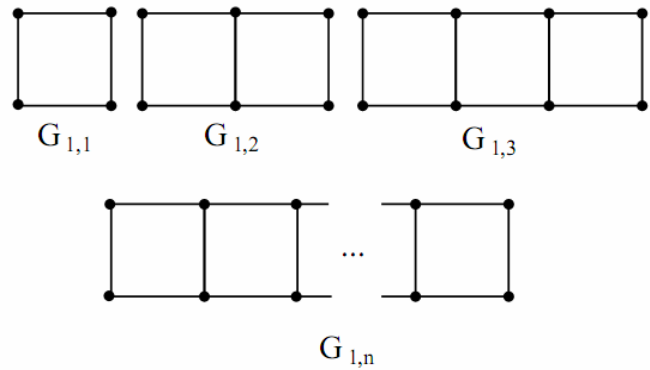

**Figure 5:** An infinite family of 2-connected graphs with all vertices inessential and four vertices of degree 2.

**Table 3.** Derivation of 2-connected graphs with *n* = 3,4,5,6 vertices from the 2-connected graphs with *n* = 3,4,5,6 vertices. Note that the name of each graph is of the form $n_k$ where n is the number of vertices and *k* is the number in the ordered list in the Appendix of [Hoover W G, De Rocco A G 1961].

| Name | Graph | → | Graph | Name | Name | Graph | → | Graph | Name |
|---|---|---|---|---|---|---|---|---|---|
| $2_1$ | | → | | $3_1$ | $4_3$ | | → | | $5_{10}$ |
| $3_1$ | | → | | $4_1$ | $5_1$ | | → | | $6_1$ |
| $3_1$ | | → | | $4_2$ | $5_2$ | | → | | $6_2$ |
| $3_1$ | | → | | $4_3$ | $5_2$ | | → | | $6_3$ |
| $4_1$ | | → | | $5_1$ | $5_2$ | | → | | $6_4$ |
| $4_2$ | | → | | $5_2$ | $5_4$ | | → | | $6_5$ |
| $4_2$ | | → | | $5_3$ | $5_4$ | | → | | $6_6$ |
| $4_2$ | | → | | $5_4$ | $5_5$ | | → | | $6_7$ |
| $4_2$ | | → | | $5_5$ | $5_5$ | | → | | $6_8$ |
| $4_3$ | | → | | $5_6$ | $5_2$ | | → | | $6_9$ |
| $4_3$ | | → | | $5_7$ | $5_5$ | | → | | $6_{10}$ |
| $4_2$ | | → | | $5_8$ | $5_6$ | | → | | $6_{11}$ |
| $4_3$ | | → | | $5_9$ | $5_6$ | | → | | $6_{12}$ |

| Name | Graph | → | Graph | Name | Name | Graph | → | Graph | Name |
|---|---|---|---|---|---|---|---|---|---|
| $5_6$ | | → | | $6_{13}$ | $5_8$ | | → | | $6_{25}$ |
| $5_4$ | | → | | $6_{14}$ | $5_2$ | | → | | $6_{26}$ |
| $5_4$ | | → | | $6_{15}$ | $5_3$ | | → | | $6_{27}$ |
| $5_5$ | | → | | $6_{16}$ | $5_7$ | | → | | $6_{28}$ |
| $5_5$ | | → | | $6_{17}$ | $5_7$ | | → | | $6_{29}$ |
| $5_5$ | | → | | $6_{18}$ | $5_7$ | | → | | $6_{30}$ |
| $5_7$ | | → | | $6_{19}$ | $5_7$ | | → | | $6_{31}$ |
| $5_7$ | | → | | $6_{20}$ | $5_8$ | | → | | $6_{32}$ |
| $5_7$ | | → | | $6_{21}$ | $5_8$ | | → | | $6_{33}$ |
| $5_7$ | | → | | $6_{22}$ | $5_9$ | | → | | $6_{34}$ |
| $5_6$ | | → | | $6_{23}$ | $5_9$ | | → | | $6_{35}$ |
| $5_8$ | | → | | $6_{24}$ | $5_5$ | | → | | $6_{36}$ |

| Name | Graph | → | Graph | Name | Name | Graph | → | Graph | Name |
|---|---|---|---|---|---|---|---|---|---|
| $5_6$ | | → | | $6_{37}$ | $5_8$ | | → | | $6_{46}$ |
| $5_6$ | | ↔ | | $6_{38}$ | $5_7$ | | → | | $6_{47}$ |
| $5_5$ | | | | | $5_{10}$ | | → | | $6_{48}$ |
| $5_6$ | | ↔ | | $6_{39}$ | $5_9$ | | → | | $6_{49}$ |
| $5_8$ | | | | | $5_9$ | | → | | $6_{50}$ |
| $5_9$ | | → | | $6_{40}$ | $5_9$ | | → | | $6_{51}$ |
| $5_9$ | | → | | $6_{41}$ | $5_8$ | | → | | $6_{52}$ |
| $5_7$ | | → | | $6_{42}$ | $5_{10}$ | | → | | $6_{53}$ |
| $5_{10}$ | | → | | $6_{43}$ | $5_9$ | | → | | $6_{54}$ |
| $5_7$ | | ↔ | | $6_{44}$ | $5_{10}$ | | → | | $6_{55}$ |
| $5_8$ | | | | | $5_{10}$ | | → | | $6_{56}$ |
| $5_8$ | | → | | $6_{45}$ | | | | | |